\documentclass{article}

\def\lb{\label}
\def\tr{{\mathrm{tr}}}
\def\dim{{\mathrm{dim}}}

\def\be{\begin{equation}}
\def\ee{\end{equation}}
\def\ba{\begin{eqnarray}}
\def\ea{\end{eqnarray}}
\def\s{\sigma}
\def\qed{\rule{5pt}{5pt}}
\def\suml{\sum\limits}

\newcommand{\vv}[2]{V_{(#1, {#2})}}
\newcommand{\rhop}[2]{\rho_{(#1, {#2})}}
\newcommand{\hh}[2]{{{\cal H}_{(#1, {#2})}}}

\begin{document}
\title{Chain models on Hecke algebra for corner type representations.}
\author{A.P.\,Isaev${}^{a}$,  O.V.\,Ogievetsky${}^{b}$  and A.F. Os'kin${}^{a}$}
\date{\empty}
\maketitle

\begin{center}
${}^{a}$ Bogoliubov Laboratory of Theoretical Physics,
Joint Institute for Nuclear Research, \\
Dubna, Moscow region 141980, Russia \\
E-mails: isaevap@theor.jinr.ru,  aoskin@theor.jinr.ru
\\
\vspace{0.3cm}
${}^{b}$ Center of Theoretical Physics\footnote{Unit\'e Mixte de Recherche
(UMR 6207) du CNRS et des Universit\'es Aix--Marseille I,
Aix--Marseille II et du Sud Toulon -- Var; laboratoire affili\'e \`a la FRUMAM (FR 2291)}, Luminy,
13288 Marseille, France \\
and P. N. Lebedev Physical Institute, Theoretical Department,
Leninsky pr. 53, 117924 Moscow, Russia \\
E-mail: oleg@cpt.univ-mrs.fr
\end{center}

\vspace{2cm}

\noindent
{\bf Abstract.}
We consider the integrable open chain models formulated in terms of generators of the Hecke
algebra. The spectrum of the Hamiltonians for the open Hecke chains of finite size with free
boundary conditions is deduced for special (corner type) irreducible representations of the Hecke
algebra.

\section{Introduction}
The A-type Hecke algebra $H_{n+1}$ is generated by the elements $\s_i$
$(i = 1,\dots, n)$ subject to the relations:
\begin{eqnarray}
\s_i \s_{i+1} \s_i &=& \s_{i+1} \s_i \s_{i+1},\label {braid}\\[-.5em]
\s_i \s_j &=& \s_i \s_j,\qquad\qquad\qquad |i - j| > 1,\label{local}\\[-.5em]
\s_i^2 &=& (q - q^{-1}) \s_i + 1,\label {hecke}
\end{eqnarray}
where $q \in {\bf C} \backslash \{0\}$ is a parameter. In this paper we consider a restricted
class of representations of $H_{n+1}$ corresponding to so called corner Young diagrams
(corner representations) related to $U_q su(1|1)$ models \cite{Ritten}. In these representations
we calculate the spectrum of the Hamiltonian of the open Hecke chain models  \cite{IsOg1, Is5},
i.e. the spectrum of the following element of the Hecke algebra
\begin{equation}
\hat{\cal H} = \sum_{i = 1}^n \s_i,\quad \s_i \in H_{n+1}. \lb{hamiltonian}
\end{equation}

The paper is organized as follows. In Section 2 we describe the corner type
representations of the Hecke algebra. In Section 3 we calculate the spectrum of the
Hamiltonian for the corner type representations for Young diagrams with two rows. In Section 4
we establish a relation between the corner type representation with $l$ rows and the $l$-th wedge
power of the corner type representation with two rows. Using this and the results of Section 3,
we calculate the spectrum of the Hamiltonian (\ref{hamiltonian}) in the general corner type
diagram.

\section{Representation for the corner type \\ diagrams $\{k+1, 1^l\}$}
The representation theory of the Hecke algebras is a well developed subject, see
\cite{Jon1, IsOgH, OgPya, Mur} and references therein. For generic $q$, it is known that
irreducible representations $\rho_{\lambda}$ of the Hecke algebra $H_{n+1}$ are labeled by the Young diagrams $\lambda$ with $n+1$ boxes and basis elements in the representation
space of $\rho_{\lambda}$ can be indexed by the standard Young tableaux of the shape $\lambda$.
The standard Young tableau for the corner diagram $\{k+1, 1^l\}$ is:

\begin{equation}\!\!\!\!\!\!\!\!\!\!\!\!\!\!\!\!\!\!\!\!\!\!\!\!\!\!\!\!\!\!\!\!\!\!\!\!\!\!\!\!
\!\!\!\!\!\!\!\!\!\!\!\!\!\!\!\!\!\!\!\!\!\!\!\!\!\!\!\!
\!\!\!\!\!\!\!\!\!\!\!\!\!\!\!\!\!\!\!\!\!\!\!\!\!\!\!\!
\begin{array}{c}
\unitlength=3mm
\begin{picture}(10,8.5)
\put(11,4){
$
\begin{array}{|c|c|c|c|c|}
\hline
1 & j_1 & j_2 & \dots & j_k \\
\hline
i_1 & \multicolumn{1}{c}{}\\
\cline{1-1}
i_2 & \multicolumn{1}{c}{}\\
\cline{1-1}
\vdots & \multicolumn{1}{c}{}\\
\cline{1-1}
i_l & \multicolumn{1}{c}{}\\
\cline{1-1}
\end{array}
$}
\put(10.5, 3.5){$\left\{\phantom{\begin{array}{c}
                                                                                                                A\\B\\C\\D
                                                                                                                \end{array}
                                                                                                                    } \right.$
}
\put(9.5, 3.5){$I$}
\end{picture}
\end{array}\lb{setidef}\end{equation}

\noindent
Here $\{ i_1,\dots ,i_l,j_1,\dots ,j_k\}$ is a $(l,k)$-shuffle of $\{ 2,3,\dots ,k+l,k+l+1\}$.
The standard Young tableau is thus determined by the set $I=\{ i_1,\dots ,i_l\}$.

Denote the corresponding representation of the Hecke algebra by $\rhop{k}{l}$ and its
space by $\vv{k}{l}$. The action of the generators $\s_p$, $1\leq p \leq k+l$, is
\begin{eqnarray}
\rhop{k}{l}(\s_p) v_I = q v_I & \mathrm{if\ } p, p + 1 \notin I,\
\lb{irrep1}\\ [0.15cm]
\rhop{k}{l}(\s_p) v_I = \frac{-q^{-p}}{(p)_q} v_I + \frac{(p - 1)_q}{(p)_q} v_{s_p I} &
\mathrm{if\ } p \notin I, p+1 \in I\ ,\lb{irrep2}\\ [0.15cm]
\rhop{k}{l}(\s_p) v_I = \frac{q^{p}}{(p)_q} v_I + \frac{(p + 1)_q}{(p)_q} v_{s_p I} &
\mathrm{if\ } p \in I, p+1 \notin I\ ,\lb{irrep3}\\ [0.15cm]
\rhop{k}{l}(\s_p) v_I = -q^{-1} v_I& \mathrm{if\ } p, p+1 \in I\ ,\lb{irrep4}
\end{eqnarray}
where $(p)_q = \frac{q^p - q^{-p}}{q - q^{-1}}$ and $v_{s_p I}$ is the basis vector corresponding
to the Young tableau with the numbers $p$ and $p+1$ interchanged.
In the matrix form ($e_{I,J} \!\in\! {\rm End}(V_{k,l})$ are the matrix units, $e_{I,J}e_{K,L}\!
=\! \delta_{JK}e_{I,L}$):
\begin{equation}
\lb{matsigm}
\begin{array}{c}
\rhop{k}{l}(\sigma_p) = q \sum\limits_{I:p,p+1 \notin I} \, e_{I,I}
+ \sum\limits_{
I:p\notin I,p+1 \in I
} \left(-\frac{q^{-p}}{p_q} \, e_{I,I}  +
\frac{(p-1)_q}{p_q} \, e_{s_p I,I} \right) \\[1em]
+ \sum\limits_{I:p \in I,p+1 \notin I} \left(\frac{q^p}{p_q} \, e_{I,I}  +
\frac{(p+1)_q}{p_q} \, e_{s_p I,I} \right)
- q^{-1} \sum\limits_{I:p,p+1 \in I} \, e_{I,I} \; .
\end{array}
\end{equation}

{\bf Proposition 2.1}
\lb{lemdim}
{\it For the Hamiltonian $\hat{\cal H} \equiv \sum\limits_{i = 1}^{k+l} \s_i\in H_{k+l+1}$ we have
\begin{equation}
\tr_{\vv{k}{l}}(\rhop{k}{l}(\hat{\cal H}))=(qk -q^{-1}l)\ \dim \vv{k}{l}\ .\lb{lemformulaitself}
\end{equation}
}

\noindent
{\bf Proof.} The dimension of the space $\vv{k}{l}$ is the number of $l$-element subsets
of $\{ 2,3,\dots ,k+l,k+l+1\}$,
$\dim \vv{k}{l} 
=\left(\!\!\!\begin{array}{c}k+l\\ l\end{array}\!\!\!\right)\ .$

By eqs.(\ref{irrep1})-(\ref{irrep4}), the action of the first generator $\s_1$ is diagonal and
\begin{equation}
\begin{array}{c}
\tr(\rhop{k}{l}(\s_1)) = q N_1 - q^{-1} N_2\lb{lemtraceformula},\,
N_1=\left(\!\!\!\begin{array}{c}k+l-1\\ l\end{array}\!\!\!\right),
N_2=\left(\!\!\!\begin{array}{c}k+l-1\\ l-1\end{array}\!\!\!\right) .
\end{array}
\end{equation}
Here $N_1$ (resp., $N_2$) is the number of sets $I$ with $2 \notin I$ (resp., $2 \in I$).

Since $X\s_iX^{-1}=\s_{i+1}\ \forall\, i \in [1, k+l - 1]$, where $X = \s_1\s_2\dots \s_{k+l}$,
the elements $\s_i,\ i > 1$, are conjugate to $\s_1$. Thus 
$\tr(\rhop{k}{l}(\s_i)) = \tr(\rhop{k}{l}(\s_1))= q N_1 - q^{-1} N_2$
and (\ref{lemformulaitself}) follows. \hfill $\qed$

It turns out that it is more convenient to work with the traceless 
Hamiltonian
\begin{equation}
\label{hamren}
\hh{k}{l}(q) := \rhop{k}{l} \left(  \sum_{i=1}^{k+l} \sigma_i \right) -
(q \, k - q^{-1} \, l){\bf 1}.
\end{equation}

\section{Spectrum of $\hh{k}{1}(q)$}\lb{spech1}
Consider the Hecke algebra $H_{k+2}$ and its representation $\rhop{k}{1}$ for the corner
Young diagram $\{k+1, 1\}$ with only two rows (i.e. $l = 1$). The dimension of this
representation  is $k+1$. As we shall see in the sequel, the Hamiltonian $\hh{k}{1}(q)$ in this
representation is a building block for the construction of the Hamiltonians $\hh{k}{l}(q)$,
corresponding to all corner diagrams $\{k+1, 1^l\}$. In the representation
$\rhop{k}{1}$ the set $I$ (\ref{setidef}) consists of only one number, $I = \{i\},\ i \in
\{2,\dots, k+2\}$, and we use the notation $v_i = v_I$ for basis vectors and $e_{i,j}$ for
matrix units.

\vskip .1cm
{\bf Proposition 3.1} {\it \lb{h1diag}
In the basis $\{ v_i \}$, the Hamiltonian (\ref{hamren}) reads 
\begin{equation}\lb{lre1}\!\!\!\!
\hh{k}{1}(q)= \suml_{p = 2}^{k+1} \! \left(\!
 \textstyle\frac{(p+1)_q e_{p, p+1} +(p - 1)_q e_{p+1, p}}{(p)_q}
-  \textstyle\frac{e_{p,p}}{(p)_q (p-1)_q}  \! \right) +
\textstyle\frac{(k)_q e_{k+2, k+2}}{(k+1)_q}  \ .
\end{equation}
}

\vskip .1cm
\noindent
{\bf Proof.} According to the general formula (\ref{matsigm}), we have
\begin{equation}\lb{sigmak1}
\rhop{k}{1}(\s_{p})  =
q {\bf 1}  + \frac{(p-1)_q}{p_q} (e_{p,p+1}-e_{p,p}) + \frac{(p+1)_q}{p_q} (e_{p+1,p}-e_{p+1,p+1}) \ \end{equation}
%where ${\bf 1} =  \sum_{i=2}^{k+2}\!\!  e_{i,i}$ and 
(here $e_{1,2}=0=e_{2,1}$).
%according to our definition of the basis $\{ v_i \}$ $(i=2,\dots,k+2)$.
Eq. (\ref{lre1}) is a straightforward consequence of (\ref{sigmak1}).\hfill$\qed$

\vskip .1cm
Let $D={\mathrm{diag}}(1,q,q^2,\dots )$. Then the operator $D\hh{k}{1}(q)D^{-1}$
possesses a finite limit ${\cal H}^{\infty}_{(k,1)}=\suml_{p = 2}^{k+1}(e_{p p+1} + e_{p+1 p})$
when $q$ tends to infinity.

For $q \in  {\bf C}^*\setminus \{ q\ |\ (k+1)_q!= 0\}$ define an upper triangular matrix $C(q)$,
\begin{equation}
C(q) = \sum_{p = 2}^{k + 2} \frac{1}{(p-1)_q} e_{p p} - \sum_{p = 2}^{k+1} \frac{1}{(p)_q} e_{p p+1} \; . \lb{isospecmatrix}
\end{equation}

{\bf Proposition 3.2}{\it \lb{propiso}
We have
\begin{equation}
\hh{k}{1}(q) C(q) = C(q) {\cal H}^{\infty}_{(k,1)}\ . \lb{isospec}
\end{equation}
}

\noindent
{\bf Proof.} A direct calculation.\hfill$\qed$

\vskip .1cm
Eq.(\ref{isospec}) demonstrates the isospectrality (the $q$-independence of the spectrum)
of the family $\hh{k}{1}(q)$. By (\ref{isospec}), $\hh{k}{1}(q)$ has the same spectrum as
${\cal H}^{\infty}_{(k, 1)}$. This spectrum is well known (${\cal H}^{\infty}_{(k, 1)}$ is the
incidence matrix of the Dynkin diagram of type $A$),
${\mathrm{Spec}}\, ({\cal H}^{\infty}_{(k,1)})=
\{2 \cos(\frac{\pi p}{k+2})\}$, $1 \leq p \leq k+1$.
We summarize the results (obtained by a different method in \cite{IsOs}).

\vskip .1cm
{\bf Theorem 3.1} {\it
The spectrum of the Hamiltonian (\ref{hamren}) for
$l=1$ is
\begin{equation}
\lb{specr1}
\mathrm{Spec} (\hh{k}{1}(q)) = \{2 \cos \frac{\pi p}{k+2}\},\ p = 1, 2, \dots, k+1.
\end{equation}
}

{\bf Remark.} Let $N=\sum_{p = 2}^{k + 2}e_{p p}-\sum_{p = 2}^{k+1}e_{p p+1}$. Then
$N^{-1}C(q)C(r)^{-1}N={\mathrm{diag}}(1,\frac{2_r}{2_q},\frac{3_r}{3_q},\dots )$. In particular,
the operators $C(q)C(r)^{-1}$ commute for different values of $q$ and $r$.

\section{Spectrum of $\hh{k}{l}(q)$}
For arbitrary $k$ and $l$ we realize the Hamiltonian $\hh{k}{l}(q)$ in terms of the
Hamiltonian $\hh{k+l-1}{1}(q)$. The best way to do this is to relate the representations
$\rhop{k}{l}$ and $\rhop{k+l-1}{1}$ of the Hecke algebra $H_{k+l+1}$.

For a vector space $V$ let $A_n$ be the antisymmetrizer in $V^{\otimes n}$ defined
by $A_n(v_1 \otimes v_2 \otimes \dots\otimes v_n) = \frac{1}{n!}\sum\limits_{s \in S_n}
 (-1)^{l(s)} v_{s(1)} \otimes v_{s(2)} \otimes \dots \otimes v_{s(n)}$
($S_n$ is the permutation group and $l(s)$ is the length of a permutation $s$). Denote
$\s_p^{(m)} = {\bf 1}^{\otimes (m-1)}\otimes\rho_{k+l-1, 1}(\s_p)\otimes
{\bf 1}^{\otimes (l-m)}\in {\mathrm{End}}(V^{\otimes l})$, where ${\bf 1}$ is the identity
matrix in $V$. Denote by $V^{\wedge l}$ the wedge power of $V$, $V^{\wedge l}
=A_lV^{\otimes l}$.

\vskip .1cm
{\bf Proposition 4.1}{\it
\lb{propofbest1}
The following identity holds
\begin{equation}
q^{1 - l}A_l(\s_p^{(1)}\s_p^{(2)}\dots\s_p^{(l)})A_l
= A_l\left(\sum_{m=1}^l \s_p^{(m)} - (l - 1)\ q\ {\bf 1}^{\otimes l}\right)A_l.  \lb{best1}
\end{equation}
}

\noindent
{\bf Proof.} The formula (\ref{best1}) is proved by induction using the $l=2$ case,
\begin{equation}\label{baseq}
q^{-1} A_2\left(\s_p^{(1)}\s_p^{(2)}\right)A_2
= A_2\left(\s_p^{(1)} + \s_p^{(2)} - q\ {\bf 1}^{\otimes 2}\right)A_2 \, ,
\end{equation}
which can be written in the form
\begin{equation}\label{baseq2}
A_2\left(\s_p^{(1)} - q \ {\bf 1} \right) \left( \s_p^{(2)}  - q \ {\bf 1} \right)A_2 = 0 \, 
\end{equation}
and directly deduced from (\ref{sigmak1}).
\hfill$\qed$

\vskip .1cm
{\bf Proposition 4.2} {\it
The set of matrices
\begin{equation}
\tilde\rho_{(k,l)}(\s_p) = q^{1-l} A_l
(\rho_{k+l-1,1}(\s_p)\otimes\dots\otimes\rho_{k+l-1,1}(\s_p)) A_l
\lb{prodrepres}
\end{equation}
defines a representation of the Hecke algebra $H_{k+l+1}$ in
$V_{k+l-1,1}^{\wedge l}$.
}

\vskip .1cm
\noindent
{\bf Proof.} The braid relations (\ref{braid}) and the locality (\ref{local})
follows from the multiplicative structure of $\tilde\rho_{(k,l)}(\s_p)$ and the fact
that $\rho_{k+l-1, 1}(\s_p)$ is a representation. The Hecke condition (\ref{hecke}) can be proved
by induction using (\ref{best1}).
\hfill$\qed$

\vskip .1cm
Proposition {\bf 4.2} can be generalized as follows.

\vskip .1cm
{\bf Proposition 4.3} {\it
Let $\rho_1$ and $\rho_2$ be representations of the Hecke algebra $H_n$ in spaces $V_1$ and $V_2$, 
respectively. Assume that an idempotent $A\in {\rm End}(V_1\otimes V_2)$, $A^2 = A$, 
commutes with $\rho_1(\s_p) \otimes \rho_2(\s_p)$ and $\rho_1(\s_p)\otimes 1 + 1\otimes \rho_2(\s_p)$ for any $p=1,2,\dots ,n-1$ and satisfies
\begin{equation}\lb{prodtosum}\textstyle
%\begin{array}{c}
%A(\rho_1(\s_p) \otimes \rho_2(\s_p)) =  (\rho_1(\s_p) \otimes \rho_2(\s_p)) A \; , \\  [0.1cm]
%A((\rho_1(\s_p)\otimes 1 + 1\otimes \rho_2(\s_p)) =
%(\rho_1(\s_p)\otimes 1 + 1\otimes \rho_2(\s_p)) A \; , \\  [0.1cm]
\left(\! (q-q^{-1} -\alpha)\rho_1(\s_p)\!\otimes\!\rho_2(\s_p)\! +\!
\rho_1(\s_p)\!\otimes\! 1\! +\! 1\!\otimes\!\rho_2(\s_p)\! +\!\frac{1-\alpha^2}{q-q^{-1}} 
1\!\otimes\! 1\!\right)\! A\! =\! 0
%\end{array}
\end{equation}
for some $\alpha\neq 0$.
Then $\rho(\s_p):= \alpha^{-1} A\, \rho_1(\s_p)\otimes \rho_2(\s_p)$ is a representation
of the Hecke algebra $H_n$ in the image $A (V_1 \otimes V_2)$ of $A$.
}

\vskip .1cm
\noindent
{\bf Proof.} A direct calculation, as in the previous Proposition. \hfill $\qed$

\vskip .1cm
The condition in (\ref{prodtosum}) factorizes as in (\ref{baseq2}) only if
$\alpha =q,-q^{-1}$. 
 
\vskip .1cm
The map $\iota :v_{i_1} \wedge v_{i_2} \wedge \dots \wedge v_{i_l}\mapsto v_I$,
$I = \{i_1, \dots ,i_{l}\}$, $i_1 < i_2 < \dots < i_l$, is an isomorphism of the
vector spaces $V_{k+l-1,1}^{\wedge l}$ and $V_{(k,l)}$ (and we use the same notation $v_I$ for
basis vectors of both spaces).
Now we identify the representation (\ref{prodrepres}) with the irreducible representation
$\rho_{(k, l)}$.

\vskip .1cm
{\bf Proposition 4.4} {\it \lb{repequiv}
The map $\iota$ intertwines the representations ${\tilde \rho_{(k,l)}}$ and $\rhop{k}{l}$.
}

\vskip .1cm
\noindent
{\bf Proof.} We directly verify eqs.(\ref{irrep1})-(\ref{irrep4})
for ${\tilde \rho_{(k,l)}}$. Note that the matrices ${\tilde \rho_{(k,l)}}(\s_p)$ can be written
in the form
\begin{equation}
\begin{array}{rcl}
\tilde\rho_{k,l}(\s_p) &=&q^{1 - l}  A_l (\rho_{k+l-1,1}(\s_p)\otimes\dots\otimes\rho_{k+l-1,1}(\s_p)) \\
&=&q^{1 - l}  (\rho_{k+l-1,1}(\s_p)\otimes\dots\otimes\rho_{k+l-1,1}(\s_p))A_l ,
\end{array}
\end{equation}
where the antisymmetrizers act from the left or from the right only.

If $p, p+1 \notin I$ then (omitting the sign of the tensor product)
\begin{equation}
\begin{array}{rcl}
\tilde\rho_{k,l}(\s_p) (v_{i_1}\wedge \dots \wedge v_{i_l}) &=&  q^{1 - l}
A_l(\s_p^{(1)} \s_p^{(2)} \dots \s_p^{(l)}) v_{i_1} v_{i_2}\dots v_{i_l} \\
&=& q A_l(v_{i_1}\dots v_{i_l}) = q v_{I},
\end{array}
\end{equation}
which proves (\ref{irrep1}).
If $p \in I, p+1 \notin I$ then
\begin{equation}\!\!\!
\begin{array}{c}
\tilde\rho_{k,l}(\s_p)(v_{i_1}\wedge \dots \wedge v_{p} \wedge \dots \wedge
v_{i_l}) = q^{1 - l}  A_l(\s^{(1)}_p\dots\s^{(l)}_p) v_{i_1}\dots v_p \dots v_{i_l} \\
= A_l \left(v_{i_1}\dots (\displaystyle\frac{q^p}{(p)_q} v_{p} + \displaystyle\frac{(p+1)_q}{(p)_q} v_{p+1}) \dots v_{i_l}\right) = \displaystyle\frac{q^p}{(p)_q} v_{I} + \displaystyle\frac{(p+1)_q}{(p)_q} v_{s_p I}\ ,
\end{array}
\end{equation}
which coincides with (\ref{irrep3}). We used that
$v_{i_1} \wedge \dots \wedge v_{p+1} \wedge \dots \wedge v_{i_l} = v_{s_p I}$
since $p+1 \notin I$. Eq.(\ref{irrep2}) can be proved in the same way. Finally, if $p, p+1 \in I$ then
\begin{equation}\!\!\!\!\!\!\!
\begin{array}{c}
\tilde\rho_{k,l}(\s_p)(v_{i_1}\wedge \dots \wedge v_{p} \wedge v_{p+1} \wedge \dots \wedge v_{i_l}) \\
=
A_l(v_{i_1}\dots (\displaystyle\frac{q^p}{(p)_q} v_p + \displaystyle\frac{(p+1)_q}{(p)_q} v_{p+1}) (-\displaystyle\frac{q^{-p}}{p_q} v_{p+1} + \displaystyle\frac{(p-1)_q}{p_q} v_p) \dots v_{i_l})  \\
=\displaystyle\frac{q^{-1}}{(p)_q^2}A_l(-v_{i_1}\dots v_pv_{p+1}\dots v_{i_l}+q^p(p-1)_qv_{i_1}\dots v_pv_p \dots v_{i_l} \\
-q^{-p}(p+1)_qv_{i_1}\dots v_{p+1}v_{p+1}\dots v_{i_l}+(p+1)_q(p-1)_qv_{i_1}\dots v_{p+1}v_{p} \dots v_{i_l}) \\
=-\displaystyle\frac{q^{-1}}{(p)_q^2}(v_{I}+(p+1)_q (p-1)_qv_{I}) = -q^{-1} v_{I},
\end{array}
\end{equation}
which coincides with (\ref{irrep4}). \hfill $\qed$

\vskip .1cm
Now one can find eigenvalues for $\hh{k}{l}(q)$ using the results of Section {\bf \ref{spech1}}.

\vskip .1cm
{\bf Theorem 4.1} {\it
The family $\hh{k}{l}(q)$ is isospectral with the spectrum
\begin{equation}\lb{maintheorem}
{\mathrm {Spec}} (\hh{k}{l}(q)) =
\left\{\suml_{i = 1}^{l} 2 \cos{\frac{\pi m_i}{k+l+1}}, \; 1 \leq m_1 < m_2 \dots < m_l \leq k+l\right\}.
\end{equation}
}

\noindent
{\bf Proof.} Due to propositions {\bf 4.1} and {\bf 4.3}, the Hamiltonian $\hh{k}{l}(q)$
equals
\begin{equation}\label{hkltosum}
\begin{array}{l}
\hh{k}{l}(q) = \left(\displaystyle\sum_{p=1}^n \tilde\rho_{k,l}(\s_p) - (q k - l q^{-1}){\bf 1}\right) A_l\\
= \left(\displaystyle\sum_{p = 1}^n \displaystyle\sum_{i=1}^{l} (\s_p^{(i)} - (l - 1) q {\bf 1}) - (k q - l q^{-1}){\bf 1}\right) A_l
 = \left( \displaystyle\sum_{i = 1}^{l} {\cal H}^{(i)} \right) A_l,
\end{array}
\end{equation}
where ${\cal H}^{(i)} = {\bf 1}^{\otimes {(i - 1)}} \otimes \hh{k+l-1}{1}(q) \otimes {\bf
1}^{\otimes (l - i)}$. The isospectrality of $\hh{k}{l}(q)$ follows from proposition {\bf 3.2}.

Let $\{\psi_m\}$, $1 \leq m \leq k+l$, be the eigenbasis of $\hh{k+l-1}{1}(q)$. By
(\ref{hkltosum}), $\{ \psi_{I}\}$, where $\psi_{I} = \psi_{m_1} \wedge \psi_{m_2} \wedge \dots
\wedge \psi_{m_l}$, $1 \leq m_1 < \dots < m_l \leq k+l$, is the eigenbasis of the Hamiltonian
$\hh{k}{l}(q)$ and (\ref{maintheorem}) follows. \hfill $\qed$

\vskip .1cm
{\bf Acknowledgement.} The first author (A.P. Isaev) was supported by
the RFBR grant No. 05-01-01086-a; the second author (O. Ogievetsky) was supported by the ANR 
project GIMP No. ANR-05-BLAN-0029-01.

\vskip .1cm
{\it Note added in proof.} After submission of this report, we became aware
of the paper G. Duchamp et al., Euler-Poincare
characteristic and polynomial representations of Iwahori-Hecke algebras, RIMS
{\bf 31} (1995), where similar results were obtained.

\end{document}